\documentclass[11pt,twoside]{article}

\usepackage{mathrsfs,amsfonts,amsmath,amssymb}
\usepackage{latexsym}
\newtheorem{satz}{Theorem}

\newtheorem{theorem}[satz]{Theorem}
\newtheorem{lemma}[satz]{Lemma}

\newtheorem{corollary}[satz]{Corollary}

\def\no{\noindent}
\def\sbeq{\subseteq}

\def\E{\mathsf {E}}

\def\F{\mathbb {F}}

\def\s{\sigma}
\def\I{{\cal I}}
\def\a{\alpha}
\def\L{{\cal L}}

\def\P{{\cal P}}

\def\a{\alpha}

\def\d{\delta}

\def\({\big (}
\def\){\big )}

\def\t{\tau}
\def\L{{\cal L}}

\def\b{\beta}
\def\ls{\leqslant}
\def\gs{\geqslant}

\begin{document}

\title{\bf On sumsets of convex sets}

\author{ By\\  \\{\sc Tomasz Schoen\footnote{The author is supported by MNSW grant N N201 543538.} ~ and ~
Ilya D.~Shkredov\footnote{This work was supported by grant RFFI NN
06-01-00383, 11-01-00759 and grant Leading ScientificSchools N
8684.2010.1.}} }
\date{}

\maketitle

\section{Introduction}
Let $A=\{a_1,\dots,a_n\},\, a_i<a_{i+1}$ be  a set of real numbers.
We say that $A$ is convex if
$$a_{i+1}-a_i>a_i-a_{i-1}$$ for every $i=2,\dots,n-1.$
Hegyv\'ari \cite{h}, answering a question of Erd{\H o}s, proved that is $A$ is convex then
$$|A+A|\gg |A|\log |A|/\log\log |A|\,.$$ This result was later improved by many authors. Konyagin \cite{k} and  Garaev \cite{g1} showed independently that
additive energy of a convex set is $\ll |A|^{5/2},$ which immediately implies that
$$|A\pm A|\gg |A|^{3/2}\,.$$
Elekes, Nathanson and Ruzsa \cite{enr} proved that if $A$ is convex then
$$|A+B|\gg |A|^{3/2}$$
for every $n$-element set $B.$
Finally, Solymosi \cite{soly} generalized the above inequality, showing that if $A$ is a set with distinct consecutive differences i.e.
$a_{i+1}-a_{i}=a_{j+1}-a_{j}$ implies $i=j$ then
\begin{equation}\label{isoly}
|A+B|\gg |A||B|^{1/2}
\end{equation}
for every set $B.$ For further  related results see \cite{g2}, \cite{g3}, \cite{gk}, \cite{ik}.
The aim of this note is to establish the following theorem.

\begin{theorem}\label{thm:main} Let $A$ be a convex set. Then
$$|A-A|\gg |A|^{8/5}\log^{-2/5} |A|
\,,
    \quad
|A+A|\gg |A|^{14/9}\log^{-2/3} |A| \,,$$ and
$$
    |A+A|^3 |A-A|^2 \log^{2} |A| \gg |A|^8 \,.
$$
\end{theorem}

In proof of Theorem 1 we  use an obvious inequality
$$|(A+A)\cap (A+A+s)|\gs |A+(A\cap (A+s))|\,,$$
which have found recently many applications \cite{kk}, \cite{sa1},
\cite{sa2}, \cite{s}, \cite{ss}. Roughly speaking, we show that if
$|A-A|$ is not much bigger than $|A|^{3/2}$ then the energy
$\E(A,A-A)$ is very large. On the other hand, we prove that this
cannot be the case for a convex set $A$.

Let us also remark that a similar result to our Theorem
\ref{thm:main} was proved independently by
 Solymosi and Szemer\'edi.

\section{Preparations}

We collect here results, which will be used in the course of the
proof of Theorem \ref{thm:main}. Denote by $\d_{A,B}(s)$ and
$\s_{A,B}(s)$ the number of representations of $s$ in the form $a-b$
and $a+b,\, a\in A,\, b\in B,$ respectively. If $A=B$ we simply
write $\d_A(s)=\d_{A,A}(s).$ Furthermore, put
$$\E(A,B)=\sum_s\d_A(s)\d_B(s)=\sum_s\d_{A,B}(s)^2=\sum_s\s_{A,B}(s)^2$$
and
$$\E_3(A)=\sum_s\d_A(s)^3\,.$$
Let $A_s=A\cap (A+s).$ Clearly, $|A_s|=\d(s).$ The following lemma was proved in \cite{ss}.

\begin{lemma}\label{le3}
    For every set $A$ we have
    \begin{equation*}
  \sum_{s} \E(A,A_s) = \E_3(A)
    \end{equation*}
\end{lemma}

\noindent Denote by $P$  and $P'$ the set of all i.e. all $s\in A-A$
such that $\d_A(s)\gs |A|^2/(2|A-A|)$ and $\d_A(s)\gs
|A|^2/(2|A+A|),$ respectively.
 Observe that
$$ \sum_{s\not\in P}|A_s|< (1/2)|A|^2\,,$$
so
\begin{equation}\label{est:P}
    \sum_{s\in P}|A_s|> (1/2)|A|^2 \,.
\end{equation}
Further
$$
    \frac{|A|^4}{|A+A|} \le \E(A,A) = \sum_{s\in A-A} |A_s|^2 \,.
$$
Thus
\begin{equation}\label{est:P'}
    \frac{|A|^4}{2|A+A|} \le \sum_{s\in P'} |A_s|^2 \,.
\end{equation}
The next is a straightforward
modification of Corollary 3.2 from \cite{ss}.

\begin{corollary}\label{corpop}
    Let $A$ be a subset of an abelian group, $P_* \subseteq A-A$ and $\sum_{s\in P_*} |A_s| = \eta |A|^2$, $\eta \in (0,1]$.
    Then
    \begin{equation*}
        \sum_{s\in P_*} |A\pm A_s| \ge \eta^2 |A|^6\E^{-1}_3(A)  \,.
    \end{equation*}
\end{corollary}
\begin{proof}
By Cauchy-Schwarz inequality we have
$$ |A| |A_s| =  \sum_x \d_{A, A_s} (x)=  \sum_x \s_{A, A_s} (x)
\ls \E(A,A_s)^{1/2}  |A\pm A_s|^{1/2}\, ,$$ so that applying once
again Cauchy--Schwarz inequality and Lemma $\ref{le3}$ we get
$$\sum_{s\in P_*}|A| |A_s|\ls \big (\sum_{s\in P_*}\E(A,A_s)\big )^{1/2} \big (\sum_{s\in
P_*}|A\pm A_s|\big )^{1/2}\ls \E_3(A)^{1/2}\big (\sum_{s\in
P_*}|A\pm A_s|\big )^{1/2}\,.~~\hfill\Box$$

\end{proof}

\bigskip

Order elements  $s\in A-A$  such that $\d_A(s_1)\gs
\d_A(s_2)\gs\dots \gs \d_A(s_t),\, t=|A-A|.$ The next lemma was
proved  in \cite{g1} and in \cite{ik}.

\begin{lemma}\label{lcon}
Let $A$ be a convex set. Then for every $r\gs 1$ we have
$$\d_A(s_r)\ll |A|/r^{1/3}\,.$$
\end{lemma}

\begin{corollary}\label{core3}
Let $A$ be a convex set. Then
$$\E_3(A)\ll |A|^3\log |A|\,.$$
\end{corollary}

The last lemma is a consequence of Szemer\'edi--Trotter theorem
\cite{sz-t}, see also \cite{ik}. We call a set $\L$ of continuous
plane curve a {\it pseudo-line system} if any two members of $\L$
share at most one point in common.

\begin{theorem}\label{thm:szemeredi-trotter} (\cite{sz-t})
Let $\P$ be a set of points and let $\L$ be a pseudo-line system.
Then
$$\I(\P,\L)=|\{(p,l)\in \P\times \L : p\in l\}|\ll |\P|^{2/3}|\L|^{2/3}+|\P|+|\L|\,.$$
\end{theorem}

\begin{lemma}\label{lpop}
Let $A$ be a convex set and $A'\sbeq A$, then for every set $B$ we
have
$$|A'+B|\gg |A'|^{3/2}|B|^{1/2}|A|^{-1/2}\,.$$
Furthermore, for every $\t\gs 1$ we have
$$|\{x\in A-B ~:~ \d_{A,B}(x)\gs \t\}|\ll \frac{|A||B|^2}{\t^3}\,.$$
\end{lemma}

\begin{proof}
Set $I=\{1,\dots,|A|\}$ and let $f$ be any continuous strictly
convex function such that $f(i)=a_i$ for every $i\in I$. Let $Q\sbeq
I$ be the set satisfying $f(Q)=A'.$ For integers $\a,\b$ put
$l_{\a,\b}=\{(q,f(q)) ~:~ q\in Q\}+(\a,\b).$ We consider  the
pseudo-line system ${\mathcal L}=\{l_{\a,\b}\}, ~\a\in I,\, \b \in
B,$ and the set  of points $\P=(Q+I)\times (A'+B).$ Thus, $\L
=|I||B|$ and $|\P|=|Q+I||A'+B|.$ Since each curve $l_{\a,\b}$
contains $|Q|$ points from $\P$ it follows by Szemer\'edi-Trotter's
theorem that
$$|Q||I||B|\ls \I (\P, \L)\ll (|\P||\L|)^{2/3}+|\P|+|\L|\,,$$
where $\I (\P, \L)$ stands for the number of incidences between $\P$
and $\L.$ Clearly, $|Q+I|\ls 2|I|=2|A|,$ so  that
$$|A'+B|\gg |A'|^{3/2}|B|^{1/2}|A|^{-1/2}\,.$$

Obviously, it is enough to prove the second assertion for $1\ll \t
\ls \min \{|A|,|B|\}$. Let $\P_\t$ be the set of points of $\P$
belonging to at least $\t$ curves from $\L$. Clearly, $\I (\P_\t,
\L)\gs \t |\P_\t|$ and by Szemer\'edi-Trotter's theorem we have
\begin{equation}\label{eq:s-t}
\t |\P_\t|\ll (|\P_\t||I|||B|)^{2/3}+|I||B|+|\P_\t|\,.
\end{equation}
We prove that $|\P_\t|\ll |I|^2|B|^2/\t^3.$ If
$(|\P_\t||I|||B|)^{2/3}\gs |I||B|$ then $|\P_\t|\ll
|I|^2|B|^2/\t^3.$ Otherwise, we have $|\P_\t|\ll |I||B|/\t\ll
|I|^2|B|^2/\t^3.$ Finally, each $x\in A+B$ with $\d_{A, B}(x)\gs \t$
gives at least $|I|$ points $p\in \P_\t$ having the same ordinate.
Therefore by (\ref{eq:s-t})
$$|\{x\in A-B: \d_{A,B}(x)\gs \t\}|\ls |\P_\t|/|I|\ll
|I||B|^2/\t^3\,,$$ which completes the proof.$\hfill\Box$
\end{proof}

\section{The proof of Theorem \ref{thm:main}}

 Put $D=A-A,\, S=A+A$ and $|D|=K|A|,\, |S|=L|A|.$
Observe that for every $s$
$$A-A_s \sbeq D\cap (D-s) \text{ ~~and~~ } A+A_s\sbeq S\cap (S+s)$$
so  $\d_D(s)\gs |A-A_s|$ and $\d_S(s)\gs |A+A_s|.$ By Corollary
\ref{corpop} and (\ref{est:P}), we have
$$\frac{|A|^6}{4E_3(A)}\ls \sum_{s\in P} |A-A_s|\ls \sum_{s\in P} \d_D(s)\,,$$
hence
$$\frac{|A|^7}{8KE_3(A)}\ls  \sum_{s} \d_A(s)\d_D(s)= \E(A,D)\,.$$
On the other hand,
$$\frac{|A|^4}{K\log |A|}\ll \E(A,D)\ll \sum_{s:\,\d_{A,D}(s)\gs \Delta} \d_{A,D}(s)^2\,,$$
where $\Delta=|A|^2/(K^2 \log |A|).$ Therefore, by Lemma \ref{lpop}
and dyadic argument, we have
$$\frac{|A|^4}{K\log |A|}\ll \sum_{j\gs 1} \Delta^2 2^{2j}\frac{K^2|A|^3}{\Delta^32^{3j}}={K^4|A|\log |A|}\,,$$
so that
$$K\gg |A|^{3/5}\log^{-2/5} |A|\,.$$

 Next we prove the second and the third inequality. Again by
 Corollary
\ref{corpop}, (\ref{est:P}) and Lemma \ref{lpop} we see that
$$\frac{|A|^7}{8KE_3(A)}\ls \sum_{s}
\d_{A,S}(s)^2\ll L^3K|A|\log |A|\,,$$
whence
$$
    |A|^8 \ll |A+A|^3 |A-A|^2 \log^{2} |A|
\,.
$$
Lemma \ref{lpop} and (\ref{est:P'}) gives
$$
    \frac{|A|^3}{L} \ll \sum_{s \in P' \,:\, |A_s| \ll L} |A_s|^2 \,.
$$
Thus, there exists $j\ge 1$, $2^j \le L^2 / |A|$ such that
$$
    \sum_{s\, :\, 2^{j-1} |A|/L < |A_s| \le 2^{j} |A|/L} |A_s| \gg \frac{|A|^2}{2^j \log |A|} \,.
$$
Hence by Corollary \ref{corpop} and Lemma \ref{lpop} we get
$$
    |A|^3 \ll L^5 2^{2j} \log^6 |A| \ll \frac{L^9}{|A|^2} \log^6
    |A|\,,
$$
so
$$
    L \gg |A|^{5/9} \log^{-2/3} |A|\,,
$$
as required. $\hfill\Box$

\bigskip

\no{Faculty of Mathematics and Computer Science,\\ Adam Mickiewicz
University,\\ Umul\-towska 87, 61-614 Pozna\'n, Poland\\} {\tt
schoen@amu.edu.pl}

\bigskip

\no{Division of Algebra and Number Theory,\\ Steklov Mathematical Institute,\\
ul. Gubkina, 8, Moscow, Russia, 119991\\} {\tt
ilya.shkredov@gmail.com}

\end{document}